\documentclass[12pt]{article} 
\newtheorem{theorem}{Theorem}
\newtheorem{proposition}[theorem]{Proposition}
\newtheorem{prop}[theorem]{Proposition}

\begin{document}
\title {Some calculus with extensive quantities:\newline wave equation}
\author{Anders Kock and Gonzalo E.\ Reyes}
\date{} \maketitle

\small ABSTRACT.  We take some first steps in providing a synthetic
theory of distributions.  In particular, we are interested in the use
of distribution theory as foundation, not just as tool, in the study
of the wave equation.


\normalsize

\section*{Introduction}
The aim of this paper is to contribute to a synthetic theory of
distributions. The sense in which we understand ``synthetic'' in this
context is that we place ourselves in a setting (category) where
everything is smooth (differentiable). Now distributions are
sometimes thought of as very non-smooth functions, like the Heaviside
function, or the Dirac delta.  We take the viewpoint, stressed by
Lawvere, that distributions are {\em extensive } quantities, where
functions are intensive ones.  It is only by a spurious comparison
with functions that distributions seem non-smooth.

\medskip
A main assumption about the category in which we work is that it is
cartesian closed, meaning that function- ``spaces'', and hence some
of the methods of functional analysis, are available.

\medskip
This viewpoint also makes it quite natural to formulate the wave
equation as an evolution equation, i.e.\ an {\em ordinary}
differential equation describing the evolution over time of any
initial distribution, so it is an ordinary differential equation with
values in the vector space of distributions.

\medskip
The main construction in the elementary theory of the wave equation
is the construction of the fundamental solution: the description of
the evolution of a point (Dirac-) distribution over time. (Other
solutions with other initial states may then by obtained by {\em
convolution} of the given initial state with the fundamental
solution; we shall not go here into this classical technique.)

\medskip
To say that {\em distributions} are extensive quantities implies that
they transform covariantly.  To say that {\em functions} are
intensive quantities implies that they transform contravariantly.
Distributions are here construed, following Schwartz, as linear
functionals on the space of (smooth) functions.  But since all
functions in the synthetic context are smooth, as well as continuous,
there is no distinction between
distributions  and Radon measures.

\medskip
So we consider a cartesian closed category ${\bf E}$ with finite
limits, in which there is given  a commutative ring object $R$, to
be thought of as the real number line.

\medskip
Already on this basis, one can define the vector space ${\cal D}'
_{c} (M)$ of distributions of compact support on $M$,  for each
object $M\in {\bf E}$, namely the object of $R$-linear maps $R^M \to
R$ (``vector space'' in this context means $R$-module).

\medskip
We shall assume that elementary differential calculus for functions
$R\to R$ is available, as in all models of SDG, cf.\  \cite{K1},
\cite{mr}, \cite{lav}, etc. We shall also assume some integral
calculus, but only in the weakest
possible sense, namely we assume

\medskip
{\bf Integration Axiom:} For every $\psi :R\to R$, there is a unique
$\Psi :R\to R$ with $\Psi ' =\psi$ and with $\Psi (0)=0$.

\medskip
Note that we do not assume any order $\leq$ on $R$, so that
``intervals'' $[a,b]\subseteq R$ do not make sense as subsets.
``Intervals'', on the contrary, will be construed as {\em
distributions}: for $a, b \in R$, $[a,b]$ denotes the distribution
$$
\psi \mapsto \int _a ^b \psi (x)\ dx
$$
The right hand side here of course means $\Psi (b)- \Psi (a)$, where
$\Psi $ is the primitive of $\psi$ given by the integration axiom.
(This weak form of integration axiom holds in some of the very simple
models of SDG, like in the topos classifying the theory of
commutative rings.)

\medskip
Finally, for the specific treatment of the wave equation, we need that the
  trigonometric functions $\cos$ and $\sin$ should be present.  We
assume that they are given as part of the data, and that they satisfy
$\cos ^2 + \sin ^2 =1$, and $\cos ' = -\sin$, $\sin ' = \cos$.  Also
as part of the data, we need specified an element $\pi \in R$ so that
$\cos \pi =-1, \cos 0 =1$.

\medskip
Except for the Taylor Series/ Formal Solutions considerations in the
end of the paper, the setting does not depend on the ``nilpotent
infinitesimals'' of SDG, but could also be, say, that of
Froelicher-Kriegl \cite{fk}, or Grothendieck's ``Smooth Topos''.

We would also like to remark that one can probably construct such
smooth toposes in which {\em no} non-trivial distribution of compact
support has a density function, or equivalently, {\em no} function
(other than $0$) gives rise to a distribution-of-compact-support; our
description of fundamental solutions to the wave equation would not
be affected. An example seems to be the topos classifying
$C^{\omega}$-algebras, where $C^{\omega}$ is the algebraic theory of
entire real- or complex- analytic functions.

\section{Generalities on distributions}
We want to apply parts of the general theory of ordinary differential
equations to some of the basic equations of mathematical physics, the
wave- and heat- equations\footnote{We do not discuss the heat
equation in the present paper; we hope to return to it and improve
the version of \cite{kr}}.

\medskip
This takes us by necessity to the realm of distributions. Not
primarily as a technique, but because of the nature of these
equations: they model evolution through time of (say) a {\em heat
distribution}. A distribution is an {\em extensive}
quantity, and does not necessarily have a density {\em function},
which is an {\em intensive} quantity; the most important of all
distributions, the point distributions (or Dirac distributions), for
instance, do not.

\medskip
As stressed by Lawvere in \cite{law},  functions are {\em intensive}
quantities, and transform {\em contravariantly}; distributions are
{\em extensive} quantities and transform {\em covariantly}.  For
functions, this is the fact that the ``space" of functions on $M$,
$R^M$ is contravariant in $M$, by elementary cartesian-closed
category theory. Similarly, the ``space" of distributions of compact
support on $M$ is
a subspace of $R^{R^M}$ (carved out by the $R$ linearity condition),
and so for similar elementary reasons is covariant in $M$.

\medskip
Let us make the formula for covariant functorality ${\cal D}' _c$
explicit.  Let $f:M\to N$ be a map.  The map ${\cal D}' _c(f):{\cal
D}' _c (M) \to {\cal D}' _c (N)$ -- which may also be denoted $f_* $
-- is described by declaring \begin{equation}<f_*  (T) , \phi > =
<T\mbox{}, \phi \circ f >, \label{functorality}\end{equation} where
$T\mbox{}$ is a distribution on $M$, and $\phi$ is a function on $N$.
The brackets denote evaluation of distributions on functions.  If we
similarly denote the value of the contravariant functor $M\mapsto
R^M$ on a map $f$ by $f^*$, the defining equation for $f_*$ goes
$<f_* (T\mbox{}),\phi > = <T\mbox{}, f^* (\phi )>$.

\medskip
We note that ${\cal D}' _c (M)$ is an $R$-linear space, and all maps
$f_* : {\cal D}'_c (M) \to {\cal D}_c' (N)$ are $R$-linear.  Also
${\cal D}'_c (M)$ is a {\em Euclidean} vector space $V$, meaning that
the basic differential calculus in available, for instance that the
basic axiom of SDG holds; we return to this in Section 2.

\medskip
For any distribution $T\mbox{}$ of compact support on $M$, one has
its {\em Total}, which is just the number $<T\mbox{},1>\in R$, where
$1$ denotes the function on $M$ with constant value 1. Since $f^* (1)
=1$ for any map $f$, it follows that $f_*$ preserves Totals.
(Alternatively, let ${\bf 1}$ denotes the terminal object (=one-point
set). Since ${\cal D}' _c ({\bf 1})\cong R$ canonically, the Total of
$T$ may also be described as $!_* (T\mbox{})$, where $! :M\to {\bf
1}$ is the unique such map. Then preservation of Totals  follows from
functorality and from uniqueness of maps into ${\bf 1}$.)

\medskip
Recall that a distribution $T$ on $M$ may be multiplied by any
function $g:M\to R$, by the recipe
\begin{equation}<g\cdot T\mbox{}, \phi > = <T\mbox{}, g\cdot \phi >.
\label{multiplying}\end{equation}

\medskip
A basic result in one-variable calculus is ``integration by
substitution''.  We present it here in pure ``distribution'' form;
note that no assumption on monotonicity or even bijectivity of the
``substitution'' $g$ is made.

\begin{prop} Given any function $g:R\to R$, and given $a,b \in R$.
Then, as distributions on $R$,
$$g_* (g' \cdot [a,b]) = [g(a),g(b)].$$\label{ibs}
\end{prop} {\bf Proof.}
Let $\psi$ be a test function, and let $\Psi$ be a primitive of it,
$\Psi' = \psi$.  So $<[g(a), g(b)], \psi> = \Psi (g(b)) - \Psi
(g(a))$.  On the other hand, by the chain rule, $\Psi \circ g$ is a
primitive of $g' \cdot (\psi \circ g)$, and so
$$
\Psi (g(a)) - \Psi (g(b)) = <[a,b], g' \cdot (\psi \circ g)> =< g'
\cdot [a,b] , \psi \circ g> = <g_* (g' \cdot [a,b]), \psi>.
$$
The external product of distributions of compact support is defined
as follows.  If $P$ is a distribution on $M$, and $Q$ a distribution
on $N$, we get a distribution $P\times Q$ on $M\times N$,  by
$$
<P\times Q, \psi >= <P , [m\mapsto <Q, \psi (m,-)>]>.
$$
In general, the external product construction $\times$ will not be
the same as the external product construction $\overline{\times}$
given by\footnote{In fact, the two external product formations
described here provide the covariant functor ${\cal D}' _{c}(-)$ with
two structures of monoidal functor ${\bf E}\to {\bf E}$, in fact,
they are
the monoidal structures that arise because ${\cal D}' _{c}(-)$ is a
strong functor on ${\bf E}$ with a monad structure, \cite{sfmm},
\cite{prsfa}.}
$$
<P\overline{\times} Q, \psi >= <Q , [n\mapsto <P, \psi (-,n)>]>.
$$
However, if $[a,b]$ and $[c,d]$ are intervals (viewed as
distributions on $R$, as described above), $[a,b]\times [c,d] =
[a,b]\overline{\times } [c,d]$, as distributions on $R^2$, by an
application of Fubini's Theorem, (which holds in the context here --
it is a consequence of equality of mixed partial dervatives).  -
Distributions arising in this way on $R^2$, we call {\em rectangles}.
The evident generalization to higher dimensions, we call {\em boxes}.
We have
$$
<[a,b]\times [c,d] , \psi > = \int _a ^b \int _c ^d \psi (x,y)\;
dy\; dx,
$$
in traditional notation. Notice that we can define {\em the boundary}
of the box $[a,b]\times [c,d]$ as the obvious distribution on $R^2$,
$$(p_c^2)_{*}[a,b] + (p_b^1)_{*}[c,d]-(p_d^2)_{*}[a,b]-(p_a^1)_{*}[c,d]$$
where $p_c^2(x)=(x,c),\; p_b^1(y)=(b,y),$ etc.

\medskip
By a {\em singular} box in an object $M$, we understand the data of a
map $\gamma :R^2 \to M$ and a box $[a,b]\times [c,d]$ in $R^2$, and
similarly for singular intervals and singular rectangles. Such a
singular box gives rise to a distribution on $M$, namely
$g_{*}([a,b]\times [c,d])$.

\medskip
By ``differential operator'' on an object $M$, we here understand
just an $R$-linear map $D: R^M \to R^M$.  If $D$ is such an operator,
and $T$ is a distribution on $M$, we define $D(T\mbox{})$ by
$$
<D(T\mbox{}), \psi >:= <T\mbox{}, D(\psi )>,
$$
and in this way, $D$ becomes a linear operator ${\cal D}' _c(M)\to
{\cal D}' _c(M).$

\medskip
In particular, if  $X$ is a vector field on $M$, one defines the
directional derivative $D_X (T\mbox{})$ of a distribution $T$ on $M$
by the formula
\begin{equation}<D_X (T\mbox{}), \psi > = <T\mbox{}, D_X (\psi )>.
\label{Lie-d}\end{equation}
This in particular applies to the vector field
$\;{\partial}/{\partial x}\;$ on $R$, and reads here
\noindent $<T\mbox{}',\psi > = <T\mbox{}, \psi ' >$ ($\psi '$ denoting the
ordinary derivative of the function $\psi$). (This is at odds with
the minus sign which is usually put in into the definition of $T'$,
but it will cause no confusion -- we
are anyway considering second order operators, where there is no discrepancy.)

\medskip
The following Proposition is an application of the covariant
functorality of the functor ${\cal D}' _c$, which  will be used in
connection with the wave equation in dimension 2. We consider the
(orthogonal) projection $p: R^3 \to R^2$ onto the $xy$-plane;
$\Delta$ denotes the Laplace operator in the relevant $R^n$, so for
$R^3$, $\Delta $ is $\partial ^2 /\partial x^2 +\partial ^2 /\partial
y^2 + \partial ^2 /\partial z^2$.

\begin{proposition} For any distribution $S$ (of compact support) on
$R^3$,
$$
p_* (\Delta (S)) = \Delta (p_* (S)).
$$
{\em (The same result holds for any orthogonal projection $p$ of
$R^n$ onto any linear subspace; the proof is virtually the same, if
one uses invariance of $\Delta$ under orthogonal transformations.)}
\label{proj}
\end{proposition}
{\bf Proof.} For any $\psi :R^2 \to R$,
$$
\Delta (p^* \psi ))= p^* (\Delta (\psi )),
$$
namely $\partial ^2 \psi /\partial x^2 + \partial ^2 \psi /\partial
y^2$.  From this, the Proposition follows purely formally.

\section{Calculus in Euclidean vector spaces}
Recall that a {\em vector space} in the present context just means an
$R$-module.  A vector space $E$ is called {\em Euclidean} if
differential and integral calculus for functions $R\to E$ is
available.  An axiomatic account is given in \cite{K1}, \cite{mr},
\cite{lav} and other places.  The coordinate vector spaces are
Euclidean, but so are also the vector spaces $R^M$, and ${\cal D}'
_{c} (M)$ for any $M$.  To describe for instance the
(``time-'')derivative $\dot{f} $ of a function $f:R\to {\cal D}' _{c}
(M)$, we put
$$
<\dot{f} (t) , \psi >= \frac{d}{dt} <f(t), \psi >.
$$
Similarly, from the integration axiom for $R$, one immediately proves
that $ {\cal D}' _c (M)$ satisfies the integration axiom, in the
sense that for any $h: R\to {\cal D}' _c (M)$, there exists a unique
$H: R\to {\cal D}' _c (M)$ satisfying $H(0)=0$ and $H'(t)= h(t) $ for
all $t$. In particular, if $h: R\to {\cal D}' _c (M)$, the
``integral'' $\int _a^b h(u)\; du$ makes sense (as $H(b)-H(a)$), and
the Fundamental Theorem of Calculus holds, almost by definition.

\medskip
As a particular case of special importance, we consider a {\em
linear} vector field on a  Euclidean $R$-module $V$. To say that the
vector field is linear is to say that its
principal-part formation $V\to V$ is a linear map, $\Gamma$, say.  We
have then the following version of a classical result.  By a {\em
formal} solution for an ordinary differential equation, we mean a
solution defined on the set $D_{\infty}$ of nilpotent elements in $R$
(these form a subgroup of $(R,+)$).

\begin{proposition} Let a linear vector field on a Euclidean vector
space $V$ be given by the linear map $\Gamma : V\to V$. Then the
unique formal solution of the corresponding differential equation,
i.e., the equation $\dot{F}(t)=\Gamma (F(t))$ with initial position
$v$, is the map $D_{\infty} \times V \to V$ given by
\begin{equation}(t,v)\mapsto e^{t\cdot \Gamma }
(v),\label{exp1}\end{equation}
where the right hand side here means the sum of the following
``series" (which has only finitely many non-vanishing terms, since
$t$ is assumed nilpotent):
$$
v + t\Gamma (v) +\frac{t^2}{2!}\Gamma ^2 (v) + \frac{t^3}{3!}\Gamma
^3 (v)+ \dots
$$
{\em (Here of course $\Gamma ^2 (v)$ means $\Gamma (\Gamma (v))$, etc.)}
\label{exp-sol}\end{proposition}
{\bf Proof.} We have to prove that $\dot{F}(t)= \Gamma (F(t))$.  We
calculate the left hand side by differentiating the series term by
term (there are only finitely many non-zero terms):
$$
\Gamma (v)
+\frac{2t}{2!}\cdot \Gamma^2 (v) + \frac{3t^2}{3!} \Gamma^3 (v) + ...  =
\Gamma (v+ t\cdot \Gamma (v) +
\frac{t^2}{2!}\cdot \Gamma ^2 (v) + ...  )
$$
using linearity of $\Gamma$.  But this is just $\Gamma$ applied to $F(t)$.

\medskip
There is an analogous result for second order differential  equations
of the form $\stackrel{\cdot \cdot}{F}(t)= \Gamma (F(t))$ (with
$\Gamma $
linear); the proof is similar and we omit it:

\begin{proposition} The formal solution of this second order
differential equation $\stackrel{\cdot \cdot}{F}= \Gamma F$, with initial
position $v$ and initial velocity $w$, is given by $$F(t) = v + t\cdot w +
\frac{t^2}{2!} \Gamma (v) + \frac{t^3}{3!} \Gamma (w) + \frac{t^4}{4!}
\Gamma ^2 (v) + \frac{t^5}{5!} \Gamma ^2
(w)+ ... .$$
\label{exp-sol2}
\end{proposition}
We shall need the following result (``change-of-variable Lemma'');
for $V=R$, it is identical to Proposition \ref{ibs}, and the proof is
in any case  the same.

\begin{prop} Given $f:R \to V$, where $V$ is a Euclidean vector
space, and given $g:R\to R$.  Then for any $a$, $b \in R$,
$$
\int _a^b f(g(x))\cdot g'(x) \; dx = \int _{g(a)}^{g(b)} f(u)\; du.
$$\label{sl}\end{prop}
Linear maps between Euclidean vector spaces preserve differentiation
and integration of functions $R\to V$; we shall explicitly need the
following particular assertion
\begin{prop}Let $F:V\to W$ be a linear map between Euclidean vector
spaces.  Then for any $f:R\to V$,
$$F(\int _a^b f(t)\; dt) =\int _a ^b F(f(t))\; dt$$.
\label{LL}\end{prop}

\section{Spheres and balls as distributions}
Let $S$ be a distribution in $R^n$; ultimately, it will be
the unit sphere, see below. We describe some families of
distributions derived from it.  Let $t \in R$ (not necessarily $t>0$
- we haven't even assumed an order relation on $R$).  We then have
the homothety ``multiplying by $t$ from $R^n$ to $ R^n$'',
which we denote $H^t$, so
$$
H^t (\underline{x}) = t\cdot \underline{x},
$$
for any $\underline{x} \in R^n$.

\medskip
We are going to use the covariant functorality of ${\cal D} ' _c$
with respect to these maps $H^t$.  Note that for any distribution $T$
on $R^n$,
\begin{equation}H^0 _* (T) = Total (T) \cdot \delta (\underline{0}),
\label{tot}\end{equation}
where $\delta (\underline{0})$ denotes the Dirac distribution at $\underline{0}
\in R^n$, given by $<\delta (\underline{0}),\psi > = \psi
(\underline{0})\;$.  We put
$$
S^t := H^t _* (S).
$$
It has the same Total as $S$, but its support\footnote{We haven't here
introduced the notion of support of a distribution, and only use the
word here for motivating the word ``diluted''.} is larger (e.g.  for
$t=2$, it is $2^{n-1}$ times as big as that of $S$). So if $S$ is the
unit sphere, $S^t$ is ``the {\em diluted sphere of radius $t$}''.  We
also want an {\em undiluted} sphere of radius $t$; we put
$$
S_t := t^{n-1}\cdot S^t.
$$
Note that in dimension 1, $S_t = S^t$.

\medskip
The {\em ball of radius 1} is made up from shells (undiluted spheres)
``of radius $u$ ($0\leq u \leq 1$)'' (heuristically !), motivating us
to put
$$
B := \int _0 ^1 S_u \; du,
$$
using integration in ${\cal D}' _c (R^n )$.  Let $t \in R$.  We put
$$
B^t := H^t_* (B).
$$
It has the same Total as $B$, but its support is larger (``if $t>1$''
- heuristically), so if $S$ is the unit sphere, $B$ is ``{\em the
diluted ball of radius} $t$'' (think of the expanding universe).  We
also want an {\em undiluted} ball of radius $t$; we put
$$
B_t := t^n \cdot B^t.
$$
We then have

\begin{prop}For all $t\in R$,
$$
B_t = \int _0 ^t S_v \; dv.
$$
\label{int}\end{prop} {\bf Proof.}
$$B_t = t^n \cdot B^t = t^n \cdot H^t _* (B)$$
$$=t^n \cdot H^t _* (\int _0 ^1 S_u \; du)$$
$$= t^n \cdot \int _0 ^1 H^t _* (S_u ) \;du$$
(by Proposition \ref{LL}) $$\;\;\;\;\;\;\;\;\;\;\;\;\;= t^n \int _0 ^1 H^t _*
(u^{n-1} \cdot
H^u _* (S))\; du$$
$$\;\;\;\;\;\;\;= t^n \cdot \int _0 ^1 u^{n-1} H^{t\cdot u}_* (S)\; du$$
$$\;\;\;\;\;\;\;\;\;\;\;\;\;=\int _0 ^1 (t\cdot u)^{n-1}\cdot
H^{t\cdot u}_* (S) \; t\; du$$
$$= \int _0 ^t v^{n-1} H^{v}_* (S)\; dv$$
(by change-of-variable Lemma (Proposition \ref{sl}), with $v:= t\cdot
u$), which is $\int _0 ^t S_v \; dv$, as claimed.

\medskip
We now give explicit defining formulae for $S$ in dimensions 1, 2 and
3.  These are of course standard integral formulae in disguise --
explicit integral formulae come by applying the definitions, and then
integral formulae for $S^t$, $S_t$ and $B_t$ may be derived (using
Proposition \ref{int} and related arguments) -- we give some of
these formulae below.

\medskip \noindent{\bf Dimension 1} $S:=\delta (1) + \delta (-1)$

\medskip \noindent{\bf Dimension 2} $S:= \mbox{cis}_* ([0, 2\pi ])$,
where $\mbox{cis}:R \to R^2$ is the map $\theta \mapsto (\cos \theta
, \sin \theta )$.

\medskip \noindent{\bf Dimension 3} $S:= s \cdot \mbox{sph} _*([0,2\pi ]
\times [0,\pi ])$, where $s: R^2 \to R$ is the function $(\theta ,
\phi )\mapsto \sin \phi $, and where $\mbox{sph}$ is ``the spherical
coordinates map'' $R^2 \to R^3$ given by
\begin{equation}(\theta , \phi )\mapsto (\cos \theta \sin \phi , \sin
\theta \sin
\phi, \cos \phi ).\label{sph}\end{equation}

In dimension 2, for instance, we have
$$
<S^t , \psi> = \int _0 ^{2\pi} \psi (t\cdot \cos \theta , t\cdot \sin
\theta )\; d\theta,
$$
$$
<S_t , \psi> = \int _0 ^{2\pi} t\cdot \psi (t\cdot \cos \theta ,
t\cdot \sin \theta )\; d\theta,
$$
and so by Proposition \ref{int},
$$
<B_t ,\psi > = \int _0 ^t [\int _0^{2\pi} u\cdot \psi (u\cdot \cos
\theta , u\cdot \sin \theta )\; d\theta ] \; du,
$$
which the reader may want to rearrange, using Fubini, into the
standard formula for integration in polar coordinates over the disk
of radius $t$; but note we have no assumptions like ``$t>0$''.

\medskip
Note also that $ B_0 = 0$, whereas $S^0$ and $B^0$ are constants
times the Dirac distribution at the origin $\underline{0}$ (use
(\ref{tot})).  The constants are the ``area'' of the unit sphere, or
the ``volume'' of the unit ball, in the appropriate dimension.
Explicitly,
\begin{equation}S^0 =2\cdot \delta (\underline{0}),\; 2\pi \cdot
\delta (\underline{0}),\; 4\pi \cdot \delta
(\underline{0}),\label{fourpi}\end{equation} and \begin{equation}B^0 =
2\cdot \delta (\underline{0}),\; \pi \cdot \delta (\underline{0}),\;
\frac{4\pi}{3}\cdot \delta (\underline{0})\label{fourthird}\end{equation}
in dimensions 1,2, and 3, respectively.

\medskip
We shall also have occasion to consider the distribution (of compact
support) $t\cdot S^t$ on $R^3$ as well as its projection $p_* (t\cdot
S^t)$ on the $xy$-plane (using functorality of ${\cal D}'_c$ with
respect to the projection map $p:R^3 \to R^2 $).

\medskip
We insert for reference two obvious ``change of variables" equations.
Recall that $H_t : R^n \to R^n$ is the homothetic transformation
``multiplying by $t$". We have, for any vector field ${\bf F}$ on
$R^n$ (viewed, via principal part, as a map $R^n \to R^n $):
\begin{equation}\mbox{div } ({\bf F}\circ H_t ) = t\cdot (\mbox{div
}{\bf F})\circ H_t ,
\label{littlediv}\end{equation}
and
\begin{equation}t^n \int_{B_1} \phi \circ H_t = \int _{B_t} \phi .
\label{ch-var}\end{equation}

\section{Divergence Theorem for Unit Sphere}
The Main Theorem of vector calculus is Stokes' Theorem: $\;\int
_{\partial \gamma} \omega = \int _{\gamma} d\omega$, for $\omega$ an
$(n-1)$-form, $\gamma$ a suitable $n$-dimensional figure (with
appropriate measure on it) and $\partial \gamma$ its geometric
boundary. In the synthetic context, the theorem holds at least for
any singular cubical chain  $\gamma :I^n \to M$ ($I^n$ the
$n$-dimensional coordinate cube), because the theorem may then be
reduced to the fundamental theorem of calculus, which is the only way
integration enters in the elementary synthetic context; measure
theory not being available therein.  For an account
of Stokes' Theorem in this context, see \cite{mr} p.139.  Below, we
shall apply the result not only for singular {\em cubes} as in
loc.cit., but also for singular {\em boxes}, like the usual ($\gamma
: R^2 \to R^2, [0, 2\pi ] \times [0,1])$, ``parametrizing the unit
disk $B$ by polar coordinates'',
\begin{equation}\gamma (\theta , r ) = (r\cos \theta , r \sin \theta ).
\label{polar}\end{equation}
We shall need from vector calculus the Gauss-Ostrogradsky ``Divergence Theorem"
$$
\mbox{flux of } {\bf F} \mbox{ over } \partial \gamma = \int
_{\gamma} (\mbox{divergence of } {\bf F}),
$$
with ${\bf F}$ a vector field, for the geometric ``figure" $\gamma $
= the unit ball in $R^n.$ For the case of the unit ball in $R^n$, the
reduction of the Divergence
Theorem to Stokes' Theorem is  a matter of the {\em differential}
calculus of vector fields, differential forms, inner products etc.\,
(See e.g. \cite{lang} p. 204). For the convenience of the reader, we
recall the case $n=2$.

\medskip
Given a vector field ${\bf F}(x,y) = (F(x,y),\; G(x,y))$ in $R^2$,
apply Stokes' Theorem to the differential form
$$
\omega := -G(x,y)dx + F(x,y)dy
$$
for the singular rectangle $\gamma$ given by (\ref{polar}) above.
Then, using the equational assumptions on $\cos, \sin$ and their
derivatives, we have
$$
\left \{
\begin{array}{lll}
\gamma ^* (dx)= \cos \theta dr - r \sin \theta d\theta \\
\gamma ^* (dy ) = \sin \theta dr + r \cos \theta d\theta \\
\gamma^* (dx\wedge dy) = r \; (dr \wedge d\theta)
\end{array}
\right.
$$
Since $\;d\omega = (\partial G/\partial y + \partial F /\partial x
)\; dx\wedge dy = \mbox{ div }({\bf F}) \; dx\wedge dy$, then
$$
\gamma^* (d\omega ) =\mbox{ div }({\bf F}) \; r \; (dr\wedge d\theta)
$$
On the other hand,
\begin{equation}\gamma ^* \omega = (F \; \sin \theta - G \; \cos
\theta )dr + (F \; r \; \cos \theta + G \; r \; \sin \theta )\;
d\theta ,\label{oneform}\end{equation}
(all $F$, $G$, and ${\bf F}$ to be evaluated at $(r\cos \theta , r\sin
\theta )$).  Therefore
$$
\int _{\gamma} d\omega = \int _0 ^{2 \pi}\int _0 ^1 \mbox{div}({\bf
F}) \; r \; dr\;  d\theta ;
$$
this is $\int _{B_1} \mbox{ div }({\bf F}) \;dA$. On the other hand
by Stokes' Theorem $\int _{\gamma} d\omega = \int _{\partial \gamma}
\omega$ which is a
curve integral of the 1-form (\ref{oneform}) around the boundary of
the rectangle $[0, 2\pi ] \times [0,1]$. This curve integral is a sum
of four terms corresponding to the four sides  of the rectangle. Two
of these (corresponding to the sides $\theta = 0$ and $\theta = 2\pi
$) cancel, and the term corresponding to the side where $r=0$
vanishes because of the $r$ in $r\; (dr\wedge d\theta)$, so only the
side with $r=1$ remains, and its contribution is, with the correct
orientation,
$$
\int _0 ^{2\pi} (F (\cos \theta , \sin \theta )\cos \theta + G (\cos
\theta , \sin \theta ) \sin \theta )\; d\theta = \int _{S_1} {\bf
F}\cdot {\bf n} \; ds
$$
where ${\bf n}$ is the outward unit normal of the unit circle. This
expression is the flux of ${\bf F}$ over the unit circle, which thus
equals the divergence integral calculated above.

\section{Time Derivatives of Expanding Spheres and Balls}
We now combine vector calculus with the calculus of the basic ball-
and sphere-distributions, as introduced in Section 3, to prove the
following result:

\begin{theorem} In $R^n$ (for any $n$), we have, for any $t$,
$$\frac{d}{dt} S^t = t \cdot \Delta(B^t),$$
($\Delta =$ the Laplace operator).
\label{prop1}\end{theorem}
{\bf Proof.} We consider the effect of the two expressions on an
arbitrary function $\psi$. We have
$$
< \frac{d}{dt} S^t , \psi> = <S, \frac{d}{dt}(\psi
\circ H^t )>\mbox{ by various definitions}
$$
$$
=<S, u \mapsto (\nabla \psi (H^t (u))\cdot u > \mbox{ by
differential calculus}
$$
$$
=   \mbox{ flux over } S \mbox{ of } (\nabla \psi
\circ H^t ) \mbox{ by special property of } S
$$
$$
= <B, div (\nabla \psi \circ H^t )> \mbox{ by divergence
Theorem}
$$
$$
= t <B, div (\nabla \psi )\circ H^t >\mbox{ by (\ref{littlediv}) }
$$
$$
=t <B,(\Delta \psi )\circ H^t >\mbox{ by definition of } \Delta
$$
$$
=t <B^t , \Delta (\psi )> = t <\Delta (B^t ), \psi > \mbox{ by various
definitions}
$$
from which the result follows.

\medskip
We collect some further information about $t$-derivatives of some of
the $t$-parametrized distributions considered.  From Proposition
\ref{int} and the Fundamental Theorem of Calculus, we immediately
derive
\begin{equation}\frac{d}{dt}(B_t )= S_t .
\label{twelve}\end{equation}
In dimension 1, we have
\begin{equation}\frac{d}{dt}(S_t )=\Delta ( B_t );
\label{sixteen}\end{equation}
for,
$$
\frac{d}{dt}<S_t ,\psi > = \frac{d}{dt}<\psi (t) + \psi (-t)> = \psi
' (t) - \psi ' (-t),
$$
wheras
$$
<\Delta B_t , \psi > = <B_t , \psi ''> = \int _{-t} ^t \psi '' (t) \; dt,
$$
and the result follows from the Fundamental Theorem of Calculus.  --
The equation (\ref{sixteen}) implies the following equation if $n=1$; we
shall prove that it also holds if $n\geq 2$: \begin{equation}t\cdot
\frac{d}{dt} (S_t ) = (n-1)S_t + t\cdot \Delta (B_t ).
\label{thirteen}\end {equation} For, differentiate $S_t = t^{n-1}\cdot S^t$
to get
$$\frac{d}{dt}(S_t ) = (n-1)t^{n-2}\cdot S^t + t^{n-1} \frac{d}{dt}(S^t
),$$ which by Theorem \ref{prop1} and the definition of $B_t$ in terms
of $B^t$  is $=(n-1)t^{n-2} \cdot S^t + \Delta (B_t )$.  Multiplying this
equation by $t$ and using the defining equation $S_t = t^{n-1}S^t$
gives the result.

\medskip
We we shall finally argue that
\begin{equation}t\cdot \frac{d}{dt}(B^t) = S^t - nB^t .
\label{fourteen} \end{equation}
For, differentiating the defining equation $B_t = t^n \cdot B^t$ gives
$d/dt B_t = nt^{n-1}\cdot B^t + t^n \cdot d/dt B^t$.  Now the left
hand side here is $S_t$, by (\ref{twelve}), so we conclude that
$t^{n-1}\cdot S^t = nt^{n-1}\cdot B^t + t^n \cdot d/dt B^t$.  If $t$
were invertible, we would conclude by cancelling $t^{n-1}$ in this
equation. But since the equation holds for {\em all} $t$, we may
cancel it in any case: a consequence of the integration axiom is the
Lavendhomme Cancellation Principle, which says that if $t\cdot
g(t)=0$ for {\em all} $t$, then $g(t)=0$ for all $t$, see \cite{lav}
Ch.1 Prop.\ 15. Applying this principle $n-1$ times then yields
(\ref{fourteen}).

\section{Wave equation}
Let $\Delta$ denote the Laplace operator $\sum \partial ^2 /\partial
x_i ^2$ on $R^n$. We shall consider the wave equation (WE) in $R^n$,
(for $n=1,2,3$),
\begin{equation}\frac{d ^2}{d t^2} Q (t) = \Delta Q(t)
\label{waveequation}\end{equation}
as a second order ordinary differential equation on the Euclidean
vector space ${\cal D}' _c (R^n )$ of distributions of compact
support; in other words, we are looking for functions $$Q: R \to
{\cal D}' _c (R^n)$$ so that for all $t\in R,$ $\ddot{Q}(t) = \Delta
(Q(t))$ (viewing $\Delta $ as a map $ {\cal D}' _c (R^n ) \to {\cal
D}' _c (R^n ) $.)

\medskip
Consider a function $f:R\to V$, where $V$ is a  Euclidean vector
space (we are interested in $V={\cal D}' _c (R^n )$) . Then we call
the pair of vectors in $V$
consisting of $f(0)$ and $\stackrel{\cdot}{f}(0)$  the {\em initial
state} of $f$. We can now, for each of the cases $n=1$, $n=3$, and
$n=2$ describe fundamental solutions to the wave equations. (The case
$n=2$ is less explicit, and is derived ``by projection'' from the one
in dimension 3.) By {\em fundamental solutions}, we mean solutions
whose initial state is either a constant times $(\delta
(\underline{0},0))$, or a constant times $(0, \delta (\underline{0})$.

\begin{theorem}
In dimension 1: The function $R \to {\cal D}' _c (R)$ given by
$$
t\mapsto
S^t (=S_t )
$$
is a solution of the WE; its  initial state is  $2(\delta
(\underline{0}), 0)$.  \newline \indent The function $R \to {\cal D}'
_c (R)$ given by
$$
t\mapsto B_t
$$
is a solution of the WE with initial state $2(0, \delta (\underline{0}))$.
\label{dim1}\end{theorem}
{\bf Proof.} We have $d/dt (B_t)=S_t$ by (\ref{twelve}), and $d/dt
(S_t) = \Delta (B_t)$, by (\ref{sixteen}).  This establishes the WE
for $B_t$.  Since $\Delta$ and $d/dt$ commute, it therefore follows
that WE also holds for $d/dt (B_t )= S_t$.  The initial position of
the solution $S_t$ is $S_0 =S^0 = 2\delta (\underline{0})$, by
(\ref{fourpi}), and the initial velocity  $\Delta (B_0 )$ by
(\ref{sixteen}), which is $0$ since $B_0 =0$. The initial state of
the solution $B_t$ is $B_0 =0$, and the initial velocity is $S_0$, as
we already calculated, so is $ = 2\delta (\underline{0})$.

\begin{theorem} In dimension 3: The function $R \to {\cal D}' _c (R^3 )$
given by
$$
t\mapsto  S^t + t^2 \Delta (B^t )
$$
is a solution of the WE with initial state $4\pi \delta
(\underline{0}),0)$. The function $R \to {\cal D}' _c (R^3 )$ given
by $$t\mapsto  t \cdot S^t $$ is a solution of the WE with initial
state $4\pi (0, \delta (\underline{0}))$.
\label{dim3}\end{theorem}
{\bf Proof.} We calculate first $d/dt$ of $t\cdot S^t$, using Theorem
\ref{prop1}: \begin{equation}\frac{d}{dt}( t\cdot S^t )= S^t + t^2
\cdot \Delta (B^t ),\label{b33}\end{equation}
and so by Theorem \ref{prop1} again,
$$
\frac{d^2 }{dt^2 } (t\cdot S^t ) =t\cdot \Delta (B^t ) +2 \cdot t
\cdot \Delta (B^t ) +t^2 \cdot \Delta (\frac{d}{dt} B^t )
$$
$$= 3\cdot t \cdot \Delta (B^t ) + t\cdot \Delta (t\cdot \frac{d}{dt}
B^t )
$$
$$
\;\;\;= 3\cdot t \cdot \Delta (B^t ) + t\cdot \Delta (S^t - 3B^t ),
$$
using (\ref{fourteen}), and now by linearity of $\Delta$, the terms
involving $\Delta (B^t )$ cancel, so we are left with the equation
\begin{equation}
\frac{d^2 }{dt^2 } (t\cdot S^t )=\Delta (t\cdot S^t ),
\label{b4}\end{equation}
which establishes WE for $t\cdot S^t$.

\medskip
Since $d/dt$ and $\Delta$ commute, and since $t\cdot S^t$ is a
solution, then so is its $t$-derivative (calculated in (\ref{b33})
above), i.e.\ $S^t + t^2 \cdot \Delta (B^t
)$ is a solution.   The assertions about initial position and
velocity follow from (\ref{fourpi}), (using Theorem \ref{prop1} to
calculate the initial velocity of the solution $S^t + t^2 \Delta (B^t
)$).

\medskip
Recall that we considered the orthogonal projection $p: R^3 \to R^2.$
Applying covariant functorality, we get for any distribution $Q$ on
$R^3$ of compact support a distribution $p_* (Q)$ on $R^2$, also of
compact support.

\begin{theorem}
In dimension 2: The function $R \to {\cal D}' _c (R^2 )$ given by
$$t\mapsto   p_* (S^t  + t^2 \Delta (B^t ) )$$
is a fundamental solution of the WE in dimension 2; its initial state
is  $4\pi (\delta (\underline{0}),0)$.The function $R
\to {\cal D}' _c (R^2 )$ given by
$$t\mapsto   p_* (t \cdot S^t )$$
is a fundamental solution of the WE in dimension 2; its initial state
is $4\pi (0,\delta (\underline{0}))$.  \end{theorem}
(Note: The $S^t$ and $B^t$ in the statement of the Theorem are those of
      $R^3$.)

\medskip
{\bf Proof.} The fact that the distributions in question are
solutions of the WE is immediate from the Proposition \ref{proj}
(``$p_* $ commutes with $\Delta $") and from the fact that $p_* :
{\cal D}'_c (R^3 )\to {\cal D}' _c (R^2 )$ is linear, and hence
commutes with formation of $d/dt$; also, ${\cal D}' _c (p)$ sends
Dirac distribution
at $\underline{0}\in R^3$ to Dirac distribution at $\underline{0}\in
R^2$, so the initial
values and velocities are as claimed.

\medskip
An explicit integral expression for the two fundamental solutions
here, obtained by projection, requires more assumptions, in
particular, a square root formation, as is known from classical
descriptions of the solutions in terms of ``Poisson's kernel''.  We
may express this by saying that the distributional solutions
presented exist under our weak assumptions, but that they are not presented
by functions (densities).

\medskip
We haven't touched the notion of support, but when defined (in a
context where it makes sense), the two fundamental solutions ,  $S^t$
and $t\cdot S^t$ in dimension 1 and 3 will have support only on the
geometric sphere of radius $t$ (which is of ``codimension'' $1$),
whereas the solution $p_* (t\cdot S^t )$ will have support in the
direct image in $R^2$ of $S^t \subseteq R^3$, and be of codimension
$0$.  This accounts for the Huygens Principle that in a 2-dimensional
world, sounds cannot be sharp signals, cf.\ e.g.\ \cite{str} p.\ 227.

\medskip
One might of course also derive one-dimensional fundamental solutions
by orthogonal projection along $q:R^3 \to R$. Since fundamental
solutions are unique modulo constants, we conclude that $q_* (S^t +
t^2 \Delta (B^t ))$ is proportional to the 1-dimensional $S^t$ (whose
support is a 2-point set).

\medskip
Combining Theorem \ref{dim3} with Proposition \ref{exp-sol2}, we can obtain
information about $S^t$, and other spheres and balls, for nilpotent
$t$.  As  examples, we shall prove
\begin{prop}If $t^5 =0$, then
in dimension 1,
$$B_t = 2 [t\cdot \delta (\underline{0}) + \frac{t^3}{3!} \cdot \delta
(\underline{0})''],
$$
and in dimension 3,
$$
S^t = 4\pi [t\cdot \delta (\underline{0}) + \frac{t^3}{3!} \cdot
\Delta (\delta (\underline{0}))].
$$
\end{prop}
{\bf Proof.} We prove the second assertion only. (The proof of the first one is
similar, using Theorem \ref{dim1}.)  We already observed in
(\ref{fourpi}) that, in dimension 3,
$$
S^0 = 4\pi \cdot \delta (\underline{0}).
$$
Now the two expressions above are both solutions to WE with initial
state $(0, 4\pi \delta (\underline{0}))$ -- the left hand side by
Theorem \ref{dim3}, and the right
hand side by Proposition \ref{exp-sol2}, with $\Gamma = \Delta$,
$v=0$, $ w= 4\pi \delta (\underline{0})$.

\end{document}